# A LIOUVIIE-TYPE THEOREMS FOR SOME CLASSES OF COMPLETE RIEMANNIAN ALMOST PRODUCT MANIFOLDS AND FOR SPECIAL MAPPINGS OF COMPLETE RIEMANNIAN MANIFOLDS

SERGEY STEPANOV

**Abstract.** In the present paper we prove Liouville-type theorems: non-existence theorems for some complete Riemannian almost product manifolds and submersions of complete Riemannian manifolds which generalize similar results for compact manifolds.

## 1. Introduction

S. Bochner devised an analytic technique to obtain non-existence theorems for some geometric objects on a closed (compact, boundaryless) Riemannian manifold, under some curvature assumption (see [1]). Currently, there are two different points of view about classical *Bochner technique*; the first one uses the *Green's divergence theorem*, and the second uses the *Hopf's theorem* which were obtained from the *Stokes's theorem* and *classical maximum principle* for compact Riemannian manifold, respectively. In particular, a good account of applications of the Bochner technique in differential geometry of Riemannian almost product manifolds and submersions may be found in [2]. We recall here that a Riemannian almost product manifold is a Riemannian manifold (*M, g*) equipped with two complementary orthogonal distributions. For instance, the total space of any submersion of an arbitrary Riemannian manifold onto another Riemannian manifold admits such a structure.

In the present paper we will use a *generalized Bochner technique*: our proofs will be based on generalized divergence theorems and a generalized maximal principle for complete, noncompact Riemannian manifolds (see [3]). We will prove Liouville type non--existence theorems for some complete, noncompact Riemannian almost product manifolds, conformal and projective diffeomorphisms and submersions of complete, noncom-

---

Stepanov Sergey

Department of Mathematics, Finance University under the Government of Russian Federation,

Leningradsky Prospect, 49-55, 125468 Moscow, Russian Federation

e-mail: *s.e.stepanov@mail.ru*

The author was supported by RBRF grant 16-01-00053-a (Russia).

pact Riemannian manifolds which generalize similar well known results for closed manifolds.

This paper is based on our report on the conference "Differential Geometry and its Applications" (July 11-15, 2016, Brno, Czech Republic).

## 2. Three global divergence theorems

Let $(M, g)$ an $n$-dimensional oriented Riemannian manifold $(M, g)$ with volume form $dVol_g = \sqrt{\det g}\, dx^1 \wedge ... \wedge dx^n$ for positively oriented local coordinates $x^1,...,x^n$. Then we can define the divergence $div\, X$ of the vector field $X$ via the formula $d(i_X\, dVol_g) = (div\, X) dVol_g$ where $i_X$ denotes contraction with $X$ (see [4, p. 281-283]). Furthermore, if $\omega$ is an $(n-1)$-form on $(M, g)$, then we can write $\omega = i_X\, dVol_g$ where $X = *_g \omega$ for the Hodge star operator relative to $g$. Thus when $\omega$ is an $(n-1)$-form with compact support in an orientable $n$-dimensional Riemannian manifold $(M, g)$ without boundary, the Stokes theorem $\int_M d\omega = 0$ follows the *classic Green divergence theorem*

$$\int_M (div\, X)\, dVol_g = 0$$

if the vector field $X$ has compact support in a (not necessarily oriented) Riemannian manifold $(M, g)$ (see [5, p. 11]). On the other hand, there are some $L^p(M, g)$-extensions of the classical Green divergence theorem to complete, noncompact Riemannian manifolds without boundary. Firstly, we formulate the following (see [6])

**Theorem 1**. *Let $(M, g)$ be geodesically complete Riemannian manifold and $X$ be a smooth vector field on $(M, g)$ which satisfies the conditions $\|X\| \in L^1(M,g)$ and $div\, X \in L^1(M,g)$, then $\int_M (div\, X)\, dVol_g = 0$ where $\|X\|$ is a norm of the vector field $X$ induced by the metric $g$.*

Later, L. Karp showed in [7] the generalized version of Theorem 1. Namely, he has provided the following theorem.

**Theorem 2**. *Let $(M, g)$ be a complete, noncompact Riemannian manifold and $X$ be a smooth vector field on $(M, g)$ which satisfies the condition*

$$\liminf_{r\to\infty} \frac{1}{r}\int_{B(r)/B(2r)} \|X\|\, dVol_g = 0 \qquad (1.2)$$

*for the geodesic ball B(r) of radius r with the center at some fixed $x \in M$. If div X has an integral (i.e. if either $(div\ X)^+$ or $(div\ X)^-$ is integrable) then $\int_M (div\ X)\, dVol_g = 0$.*

In particular, from the above theorem we conclude that if outside some compact set *div X* is everywhere $\geq 0$ (or $\leq 0$) then $\int_M (div\ X)\, dVol_g = 0$.

In conclusion, we formulate the third generalized Green's divergence theorem (see [8]; [9]), which can be regarded, as a consequence of the above two theorems and Yau lemma from [10].

**Theorem 3**. *Let X be a smooth vector field on a connected complete, noncompact and oriented Riemannian manifold (M, g), such that $div\ X \geq 0$ (or $div\ X \leq 0$) everywhere on (M, g). If the norm $\|X\| \in L^1(M,g)$, then $div\ X = 0$.*

The *Laplace-Beltrami operator* of any $f \in C^\infty M$ is defined as $\Delta f = div_g(grad\ f)$ where *grad f* is the unique vector field that satisfies $g(X, grad\ f) = df(X)$ for all vector fields *X* on *M*. The scalar function $f \in C^\infty M$ is said to be *harmonic* if it satisfies $\Delta f = 0$. It is well known, that if $f \in C^\infty M$ is a harmonic function on any complete Riemannian manifold satisfying $f \in L^p(M,g)$ for some $1 < p < \infty$, then *f* is constant (see [10]).

In addition, we recall that the scalar function $f \in C^\infty M$ is called *subharmonic* (resp. *superharmonic*) if $\Delta f \geq 0$ (resp. $\Delta f \leq 0$). In particular, if (*M, g*) is compact then every harmonic (subharmonic and superharmonic) function is constant by the Hopf's theorem. On the other hand, Yau has proved in [10] that any subharmonic function $f \in C^2 M$ defined on a complete, noncompact Riemannian manifold with $\int_M \|df\| dV_g < \infty$ is harmonic. Then based on this statement (or on the Theorem 3) we conclude that the following lemma is true.

**Lemma**. *If (M, g) is a complete, noncompact Riemannian manifold, then any superharmonic function $f \in C^2 M$ with gradient in $L^1(M,g)$ is harmonic.*

**Proof.** Let $(M, g)$ be a complete Riemannian manifold and $f$ be a scalar function such that $f \in C^2 M$, $\Delta f \leq 0$ and $\|\operatorname{grad} f\| \in L^1(M, g)$. If we suppose that $\varphi = -f$ then the above conditions can be written in the form $\varphi \in C^2 M$, $\Delta \varphi \geq 0$ and $\|\operatorname{grad} \varphi\| \in L^1(M, g)$. In this case, from the Yau statement we conclude that $\Delta \varphi = 0$ and hence $f = -\varphi$ is a harmonic function.

## 3. Liouville-type theorems for some complete Riemannian almost product manifolds

Let $(M, g)$ be an $n$-dimensional $(n \geq 2)$ Riemannian manifold with the Levi-Civita connection $\nabla$ and $TM = \mathcal{V} \oplus \mathcal{H}$ be an orthogonal decomposition of the tangent bundle $TM$ into *vertical* $\mathcal{V}$ and *horizontal* $\mathcal{H}$ distributions of dimensions $n - m$ and $m$, respectively. We shall use the symbols $V$ and $H$ to denote the orthogonal projections onto $\mathcal{V}$ and $\mathcal{H}$, respectively. In this case we can define a *Riemannian almost product manifold* (see [11]) as the triple $(M, g, P)$ for $P = V - H$, where $(M, g)$ is a Riemannian manifold $M$ and $P$ is a (1,1)-tensor field on $M$ satisfying $P^2 = \operatorname{id}$ and $g(P, P) = g$. In addition, the eigenspaces of $P$ corresponding to the eigenvalues $1$ and $-1$, at each point, determine two orthogonal complementary distributions $\mathcal{V}$ and $\mathcal{H}$.

The *second fundamental form* $Q_H$ and *the integrability tensor* $F_H$ of $\mathcal{H}$ are define by (see [12, p. 148])

$$Q_H = \frac{1}{2} V(\nabla_{HX} HY + \nabla_{HY} HX), \qquad F_H = \frac{1}{2} V(\nabla_{HX} HY - \nabla_{HY} HX)$$

for any smooth vector fields $X$ and $Y$ on $M$. It is well known that $Q_H$ vanishes if and only if $\mathcal{H}$ is a *totally geodesic distribution*. We recall that a distribution on a Riemannian manifold is totally geodesic if each geodesic which is tangent to it at point remains for its entire length (see [12, p. 150]). On the other hand, $F_H$ vanishes if and only if $\mathcal{H}$ is an integrable distribution. A maximal connected integral manifold of $\mathcal{H}$ is called a *leaf* of the foliation. The collection of leafs of $\mathcal{H}$ is called a *foliation* of $M$. By interchanging and we define the corresponding tensor fields $Q_V$ and $F_V$ for $\mathcal{V} = \mathcal{H}^\perp$.

We define now the *mixed scalar curvature* of $(M, g)$ as the function

$$s_{\text{mix}} = \sum_{a=1}^{m} \sum_{\alpha=m+1}^{n} \sec(E_a, E_\alpha)$$

where $\sec(E_a, E_\alpha)$ is the sectional curvature of the *mixed plane* $\pi = \text{span}\{E_a, E_\alpha\}$ for the local orthonormal frames $\{E_1,...,E_m\}$ and $\{E_{m+1},...,E_n\}$ on *TM* adapted to $\mathcal{V}$ and $\mathcal{H}$, respectively (see [11]; [13, p. 23]). It is easy to see that this expression is independent of the chosen frames. With this in hand we can now state the formula which can be found in [14] and [15]. Namely, the following formula:

$$div(\xi_V + \xi_H) =$$
$$= s_{\text{mix}} + \|Q_V\|^2 + \|Q_H\|^2 - \|F_V\|^2 - \|F_H\|^2 - \|\xi_V\|^2 + \|\xi_H\|^2 \quad (2.1)$$

where $\xi_V = \text{trace}_g Q_V$ and $\xi_H = \text{trace}_g Q_H$ are the *mean curvature vectors* of $\mathcal{V}$ and $\mathcal{H}$, respectively (see [11]; [12, p. 149]).

Assume that $\mathcal{V}$ and $\mathcal{H}$ are *totally umbilical distributions*, i.e. $Q_V = (n-m)^{-1} g(V,V) \otimes \xi_V$ and $Q_H = m^{-1} g(H,H) \otimes \xi_H$ (see [11]; [12, p. 151]). In this case the formula (2.1) can be rewrite in the form (see [15]; [16]; [17])

$$div(\xi_V + \xi_H) = s_{\text{mix}} - \|F_V\|^2 - \|F_H\|^2 - \frac{n-m-1}{n-m}\|\xi_V\|^2 + \frac{m-1}{m}\|\xi_H\|^2 \quad (2.2)$$

If in addition to the assumption above we now suppose that (*M, g*) is a connected complete and oriented Riemannian manifold without boundary and with nonpositive mixed scalar curvature $s_{\text{mix}}$, then from (3.2) we obtain the inequality $div(\xi_V + \xi_H) \leq 0$. If at the same time, $\|\xi_V + \xi_H\| \in L^1(M,g)$ then by Theorem 4 we conclude that $div(\xi_V + \xi_H) = 0$. In this case, from (3.2) we obtain $\|\xi_H\|^2 = \|\xi_V\|^2 = \|N\|^2 = 0$. It means that $\mathcal{V}$ and $\mathcal{H}$ are two integrable distributions with totally geodesic integral manifolds (totally geodesic foliations). We fix now a point $x \in M$ and let $M_1$ and $M_2$ be the maximal integral manifolds of distributions through *x*, respectively. Then by the *de Rham decomposition theorem* (see [4, p. 187]) we conclude that if (*M, g*) is a simply connected Riemannian manifold then it is isometric to the direct product $(M_1 \times M_2, g_1 \oplus g_2)$ of some Riemannian manifolds $(M_1, g_1)$ and $(M_2, g_2)$ for the Riemannian metric $g_1$ and $g_2$ which induced by *g* on $M_1$ and $M_2$, respectively. In addition, we recall that every simply connected mani-fold $M$ is orientable. Summarizing, we formulate the statement which generalizes a theorem on two orthogonal complete totally umbilical distributions on compact Riemannian manifold with non positive mixed scalar curvature that has been proved in [11]; [15]; [16] and [18].

**Theorem 4**. *Let $(M, g)$ be a complete, noncompact and simply connected Riemannian manifold with two orthogonal complementary totally umbilical distributions $\mathcal{V}$ and $\mathcal{H}$ such that their mean curvature vectors $\xi_V$ and $\xi_H$ satisfy the condition $\|\xi_V + \xi_H\| \in L^1(M, g)$. If the mixed scalar curvature $s_{\mathrm{mix}}$ of $(M, g)$ is nonpositive then $\mathcal{V}$ and $\mathcal{H}$ are integrable and $(M, g)$ is isometric to a direct product $(M_1 \times M_2, g_1 \oplus g_2)$ of some Riemannian manifolds $(M_1, g_2)$ and $(M_1, g_2)$ such that integral manifolds of $\mathcal{V}$ and $\mathcal{H}$ correspond to the canonical foliations of the product $M_1 \times M_2$.*

We consider now an $(n-1)$-dimensional totally geodesic distribution $\mathcal{H}$ on $(M, g)$. In this case the formula (2.2) can be rewrite in the form (see [15])

$$\mathrm{div}\, \xi_V = s_{\mathrm{mix}} - \|F_H\|^2 \tag{2.2}$$

where $s_{\mathrm{mix}} = \sum_{\alpha=2}^{m} \mathrm{sec}(e_1, e_\alpha) = \mathrm{Ric}(e_1, e_1)$ is the *vertical Ricci curvature* for an orthonormal frame $\{e_1, \ldots, e_n\}$ at a point $x \in M$ such that $\mathcal{H}_x = \mathrm{span}\{e_1\}$ and $\mathcal{V}_x = \mathrm{span}\{e_2, \ldots, e_n\}$. Hence, an immediate consequence of (2.2) and Theorem 4 is following

**Corollary 1**. *Let $(M, g)$ be an n-dimensional complete noncompact and simply connected Riemannian manifold with $(n-1)$-dimensional totally geodesic horizontal distribution $\mathcal{H}$. If the vertical Ricci curvature of $(M, g)$ is nonpositive and $\|\xi_V\| \in L^1(M, g)$ for the mean curvature vector of $\mathcal{V}$, then $\mathcal{H}$ is integrable and $(M, g)$ is isometric to a direct product $(M_1 \times M_2, g_1 \oplus g_2)$ of Riemannian manifolds $(M_1, g_2)$ and $(M_2, g_2)$ such that $\dim M_1 = 1$ and integral manifolds of $\mathcal{V}$ and $\mathcal{H}$ correspond to the canonical foliations of the product $M_1 \times M_2$.*

The integral formula (2.1) can be reformulated as follows (see [11])

$$2\mathrm{div}_\mathcal{V}\, \xi_H + 2\mathrm{div}_\mathcal{H}\, \xi_V = 4s_{\mathrm{mix}} + \frac{1}{2}\|\nabla P\|^2 - \|F_V\|^2 - \|F_H\|^2. \tag{2.3}$$

where $div_\mathcal{V}\, \xi_H = \sum_{a=1}^{m} g(\nabla_{E_a}\xi_H, E_a)$ and $div_\mathcal{H}\, \xi_V = \sum_{\alpha=m+1}^{n} g(\nabla_{E_\alpha}\xi_V, E_\alpha)$. If $\mathcal{V}$ and $\mathcal{H}$ are integrable distributions then (2.3) can be rewrite in the form

$$2 div_\mathcal{V}\, \xi_H + 2 div_\mathcal{H}\, \xi_V = 4 s_{mix} + \frac{1}{2}\|\nabla P\|^2. \qquad (2.4)$$

Suppose now that all integral manifolds of the vertical distribution $\mathcal{V}$ are minimal submanifolds of the Riemannian manifold $(M, g)$ and $s_{mix} \geq 0$. Then from (2.4) we obtain

$$div_\mathcal{V}\, \xi_H = 2 s_{mix} + \frac{1}{4}\|\nabla P\|^2 \geq 0. \qquad (2.5)$$

If at least one connected complete and oriented maximal integral manifold $M'$ of $\mathcal{V}$ exists. We assume that $M'$ equipped with the Riemannian metric $g'$ inherited from $(M, g)$ such that $\|\xi_H\| \in L^1(M', g')$ for the mean curvature vector $\xi_H$ of $\mathcal{H} = \mathcal{V}^\perp$ which belongs to $M'$ at each point $x \in M'$. Then by applying Theorem 3 to $\xi_H$, from (2.5) we get $div_{g'}\, \xi_H = 0$. Therefore, if all integral manifolds of $\mathcal{V}$ are connected complete and oriented minimal submanifolds of the Riemannian manifold $(M, g)$ and $\xi_H$ is a $L^1$-vector field for every of them, then $\nabla P = 0$. In this case $\mathcal{V}$ and $\mathcal{H}$ are two integrable distributions with totally geodesic integral manifolds (totally geodesic foliations) on $(M, g)$ (see [11]). If at the same time, $(M, g)$ is complete, noncompact and simply connected Riemannian manifold then by the *de Rham decomposition theorem* (see [4, p. 187]) it is isometric to the direct product $(M_1 \times M_2,\, g_1 \oplus g_2)$ of some Riemannian manifolds $(M_1, g_1)$ and $(M_2, g_2)$ for the Riemannian metric $g_1$ and $g_2$ which induced by $g$ on $M_1$ and $M_2$, respectively. Summarizing, we formulate the statement which generalizes the main theorem of [12].

**Theorem 5.** *Let $(M, g)$ be complete, noncompact and simply connected Riemannian manifold. If the following three conditions are satisfied:*

1) *$(M, g)$ admits an integrable distribution $\mathcal{V}$ such that an arbitrary integral manifold $(M', g')$ of $\mathcal{V}$ which equipped with the Riemannian metric $g'$ inherited from $(M, g)$ is a connected complete and oriented minimal submanifold of $(M, g)$;*

2) *the orthogonal complementary distribution* $\mathcal{H} = \mathcal{V}^\perp$ *is also integrable and its mean curvature vectors* $\xi_H$ *satisfies the condition* $\|\xi_H\| \in L^1(M', g')$;

3) *the mixed scalar curvature* $s_{mix} \geq 0$, *then* $(M, g)$ *is isometric to a direct product* $(M_1 \times M_2, g_1 \oplus g_2)$ *of some Riemannian manifolds* $(M_1, g_2)$ *and* $(M_1, g_2)$ *such that integral manifolds of* $\mathcal{V}$ *and* $\mathcal{H}$ *correspond to the canonical foliations of the product* $M_1 \times M_2$.

**Remark 1.** If in addition, at least one closed integral manifold $M'$ of $\mathcal{V}$ exists, then, by applying the classic Green divergence theorem to $\xi_H$, from (2.5) we get

$$\int_{M'} \left( 8 s_{mix} + \|\nabla P\|^2 \right) dVol_{g'} = 0$$

where $dVol_{g'}$ is the volume form of $(M', g')$. If $(M', g')$ is non-oriented we can consider its orientable double cover. In this case the inequality $s_{mix} > 0$ is a condition of nonexistence of two orthogonal complementary foliations one of which consists of minimal submanifolds.

### 3. Applications to the theory of projective mappings of Riemannian manifolds

We recall here the definition of *pregeodesic* and *geodesic curves*. Namely, a *pregeodesic curve* is a smooth curve $\gamma: t \in J \subset \mathbb{R} \to \gamma(t) \in M$ on a Riemannian manifold $(M, g)$, which becomes a geodesic curve after a change of parameter. Let us change the parameter along $\gamma$ so that $t$ becomes an *affine parameter*. Then $\nabla_X X = 0$ for $X = d\gamma/dt$, and $\gamma$ is called a *geodesic curve*. By analyzing of the last equation, one can conclude that either $\gamma$ is an immersion, i.e., $d\gamma/dt \neq 0$ for all $t \in J$, or $\gamma(t)$ is a point of $M$.

Let $(M, g)$ and $(\overline{M}, \overline{g})$ be Riemannian manifolds of dimension $n \geq 2$. Then a smooth map $f : (M, g) \to (\overline{M}, \overline{g})$ of Riemannian manifolds is a *projective map* if $f(\gamma)$ is a pregeodesic in $(\overline{M}, \overline{g})$ for an arbitrary pregeodesic $\gamma$ in $(M, g)$ (see [19]). In particular, if a projective map $f : (M, g) \to (\overline{M}, \overline{g})$ is called *totally geodesic* if it maps linearly parametrized geodesics of $(M, g)$ to linearly parametrized geodesics of $(\overline{M}, \overline{g})$. An equiva-lent definition is that *f* is connection-preserving, or *affine*. The global structure of these maps is investigates in the paper [20].

For a projective diffeomorphism $f:(M,g)\to(\overline{M},\overline{g})$ we have (see [21, p. 135])

$$\overline{Ric} = Ric + (n-1)(\nabla d\psi - d\psi \otimes d\psi) \tag{3.1}$$

where *Ric* and $\overline{Ric}$ denote the Ricci tensors of (*M*, *g*) and $(\overline{M},\overline{g})$, respectively, and

$$\psi = \frac{1}{2(n+1)}\log\left(\frac{\det \overline{g}}{\det g}\right) + C \tag{3.2}$$

for some constant *C*. Now we can formulate the following

**Theorem 6**. *Let (M, g) be a connected complete, noncompact Riemannian manifold and $f:(M,g)\to(\overline{M},\overline{g})$ be a projective diffeomorphism onto another Riemannian manifold $(\overline{M},\overline{g})$ such that $trace_g \overline{Ric} \geq s$ for the Ricci tensor $\overline{Ric}$ of $(\overline{M},\overline{g})$ and the scalar curvature s of (M, g). If the gradient of the function $\log\left(\frac{\det \overline{g}}{\det g}\right)$ has integrable norm on (M, g) then f is affine map.*

**Proof.** We conclude immediately from (3.1) that

$$\Delta\psi = \frac{1}{n-1}\left(trace_g \overline{Ric} - s\right) + \|grad\,\psi\|^2 \tag{3.3}$$

Let $trace_g \overline{Ric} \geq s$ then (3.3) shows $\Delta\psi \geq 0$. If (*M*, *g*) is a complete, noncompact Riemannian manifold and $\|grad\,\psi\| \in L^1(M,g)$ then by the Yau statement (see [10, p. 660]) we conclude that $\Delta\psi = 0$ and $\psi$ must be harmonic on (*M*, *g*). At the same time, we see from (3.3) that $\psi$ is constant. Then according to the formula (40.8) from [21, p. 133] we obtain $\nabla \overline{g} = 0$. Hence by [20], *f* is affine map.

Let (*M*, *g*) and $(\overline{M},\overline{g})$ be Riemannian manifolds of dimension *n* and *m* such that *n* > *m*. A surjective map $f:(M,g)\to(\overline{M},\overline{g})$ is a *submersion* if it has maximal rank *m* at any point *x* of *M*, that is, each differential map $f_{*x}$ of *f* is surjective, hence, for $y \in \overline{M}$. In this case, $f^{-1}(y)$ for an arbitrary $y \in \overline{M}$ is an (*n* – *m*)-dimensional closed submanifold *M'* of (*M*, *g*) (see [22, p.11]). We call the submanifolds $f^{-1}(y)$ *fibers*.

Putting $\mathcal{V}_x = Ker(f_*)_x$, for any $x \in M$, we obtain an integrable vertical distribution $\mathcal{V}$ which corresponds to the foliation of $M$ determined by the fibres of $f$, since each $\mathcal{V}_x = T_x f^{-1}(y)$ coincides with tangent space of $f^{-1}(y)$ at $x$ for $f(x) = y$.

Let $\mathcal{H}$ be the complementary distribution of $\mathcal{V}$ determined by the Riemannian metric $g$, i.e. $\mathcal{H}_x = \mathcal{V}_x^\perp$ at each $x \in M$. So, at any $x \in M$, one has the orthogonal decomposition $T_x(M) = \mathcal{V}_x \oplus \mathcal{H}_x$ where $\mathcal{H}_x$ is called the *horizontal space* at $x$. Thus we have defined a Riemannian almost product structure on $(M, g)$.

Consider now an $n$-dimensional simple connected complete Riemannian manifold $(M, g)$, and suppose that a projective submersion $f : (M, g) \to (\overline{M}, \overline{g})$ onto an $m$-dimensional ($m < n$) Riemannian manifold $(\overline{M}, \overline{g})$ exists. Then each pregeodesic line $\gamma \subset M$ which is an integral curve of the distribution $Ker f_*$ is mapped into a point $f(\gamma)$ in $\overline{M}$. Note that this fact does not contradict the definition of projective submersion.

In addition, we have proved in [23] and [24] that $(M, g)$ is isometric to a *twisted product* $\left(M_1 \times M_2, g_1 + e^{2\alpha_2} g_2\right)$ of some Riemannian manifolds $(M_1, g_2)$ and $(M_2, g_2)$, and for smooth function $\alpha_2 : M_1 \times M_2 \to \mathbb{R}$ such that all fibres of submersion and their orthogonal complements correspond to the canonical foliations of $M_1 \times M_2$ (see [23] and [24]). In this case, the following corollary of Theorem 4 is true.

**Corollary 2**. *Let $(M, g)$ be an n-dimensional complete, noncompact and simply connected Riemannian manifold and $f : (M, g) \to (\overline{M}, \overline{g})$ be a projective submersion onto another m-dimensional ($m < n$) Riemannian manifold $(\overline{M}, \overline{g})$ such that the mean curvature vector $\xi_H$ of the horizontal distribution $(Ker f_*)^\perp$ satisfies the condition $\|\xi_H\| \in L^1(M, g)$. If the mixed scalar curvature $s_{mix}$ of $(M, g)$ is nonpositive then $(Ker f_*)^\perp$ is integrable and $(M, g)$ is isometric to a direct product $(M_1 \times M_2, g_1 \oplus g_2)$ of some Riemannian manifolds $(M_1, g_2)$ and $(M_1, g_2)$ such that integral manifolds of $Ker f_*$ and $(Ker f_*)^\perp$ correspond to the canonical foliations of the product $M_1 \times M_2$.*

Moreover, we have proved in [24] that if a simple connected complete $n$-dimensional Riemannian manifold $(M, g)$ has a nonnegative sectional curvature and admits a projec-tive submersion onto another *m*-dimensional (*m* < *n*) Riemannian manifold $(\overline{M}, \overline{g})$, then (*M, g*) is isometric to a direct product $(M_1 \times M_2, g_1 \oplus g_2)$ of some Riemannian manifolds $(M_1, g_2)$ and $(M_2, g_2)$ such that the integral manifolds of $Ker\, f_*$ and $(Ker\, f_*)^\perp$ correspond to the canonical foliations of the product $M_1 \times M_2$. We can formulate now a statement which will supplement this theorem. The statement is a corollary of Theorem 1 and Theorem 6.

**Corollary 3**. *Let (M, g) be an n-dimensional complete, noncompact and simply connected Riemannian manifold and $f: (M, g) \to (\overline{M}, \overline{g})$ be a projective submersion onto another m-dimensional (m < n) Riemannian manifold $(\overline{M}, \overline{g})$ with connected fibres. If the mixed scalar curvature $s_{mix} \geq 0$ then (M, g) is isometric to a direct product $(M_1 \times M_2, g_1 \oplus g_2)$ of some Riemannian manifolds $(M_1, g_2)$ and $(M_2, g_2)$ such that the integral manifolds of $Ker\, f_*$ and $(Ker\, f_*)^\perp$ correspond to the canonical foliations of the product $M_1 \times M_2$.*

**Proof.** Let (*M, g*) be an *n*-dimensional complete, noncompact and simply connected Riemannian manifold and $f: (M, g) \to (\overline{M}, \overline{g})$ be a projective submersion onto another *m*-dimensional (*m* < *n*) Riemannian manifold $(\overline{M}, \overline{g})$ with connected fibres. It follows from the above, the fibre $f^{-1}(y)$ for an arbitrary $y \in \overline{M}$ is an a (*n* – *m*)-dimensional closed connected submanifold $M'$ of (*M, g*) equipped with the Riemannian metric $g'$ inherited from (*M, g*). The mean curvature vector $\xi_H$ of $\mathcal{H} = (Ker\, f_*)^\perp$ belongs to $T_x M'$ at each point $x \in M'$ then, by applying the classic Green divergence theorem $\int_{M'} (div_\mathcal{V}\, \xi_H)\, d\,Vol_{g'} = 0$ to $\xi_H$, from (2.4) we get the following equation

$$\int_{M'} \left(2 s_{mix} + 1/4 \|\nabla P\|^2\right) d\,Vol_{g'} = 0 \tag{3.4}$$

where $dVol_{g'}$ is the volume form of $(M', g')$. If $(M', g')$ is non-oriented we can consider its orientable double cover. If the mixed scalar curvature $s_{mix} \geq 0$ then from (3.4) we obtain that $\nabla P = 0$ at each point of $M'$. At the same time, we recall that $M'$ is an arbitrary fibre of the projective submersion $f: (M, g) \to (\overline{M}, \overline{g})$. Therefore, $Ker\, f_*$ and $(Ker\, f_*)^\perp$ are two integrable distributions with totally geodesic integral manifolds (totally geodesic foliations) on the complete, noncompact and simply connected Riemannian manifold $(M, g)$. Then by the well known de Rham decomposition theorem it is isometric to the direct product $(M_1 \times M_2, g_1 \oplus g_2)$ of some Riemannian manifolds $(M_1, g_1)$ and $(M_2, g_2)$ for the Riemannian metric $g_1$ and $g_2$ which induced by $g$ on $M_1$ and $M_2$, respectively. The proof of our corollary is complete.

Using the equality (3.4) once again we get the following

**Corollary 4.** *Let $(M, g)$ be an n-dimensional Riemannian manifold and $f:(M,g) \to (\overline{M},\overline{g})$ be a submersion onto another m-dimensional ($m < n$) Riemannian manifold $(\overline{M},\overline{g})$ with connected fibres. If the mixed scalar curvature $s_{mix} > 0$ then $f:(M,g) \to (\overline{M},\overline{g})$ is not a projective submersion.*

**Proof.** For the case $s_{mix} > 0$ we can rewrite (3.4) in the form

$$\int_{M'} \left(2 s_{mix} + 1/4 \|\nabla P\|^2 \right) d Vol_{g'} > 0.$$

The contradiction just obtained with the classic Green divergence theorem completes the proof of our corollary.

### 4. Applications to the theory of conformal mappings of Riemannian manifolds

Let $(M, g)$ and $(\overline{M},\overline{g})$ be Riemannian manifolds of dimension $n \geq 2$. Then a diffeomorphism $f:(M,g) \to (\overline{M},\overline{g})$ is called *conformal* if it preserves angles between any pair curves. In this case, $\overline{g} = e^{2\sigma} g$ for some scalar function $\sigma$. In particular, if the function $\sigma$ is constant then $f$ is a *homothetic mapping*.

If $\sigma \in C^2 M$ then (see [21, p. 90])

$$e^{2\sigma} \overline{s} = s - 2(n-1) \Delta \sigma - (n-1)(n-2) \|\text{grad }\sigma\|^2 \qquad (4.1)$$

where $\overline{s}$ denote the scalar curvature $(\overline{M},\overline{g})$. Now we can formulate the following

**Theorem 6.** *Let $(M, g)$ be an n-dimensional ($n \geq 3$) complete, noncompact Riemannian manifold and $f:(M,g) \to (\overline{M},\overline{g})$ be a conformal diffeomorphism onto another Riemannian manifold $(\overline{M},\overline{g})$ such that $\overline{g} = e^{2\sigma} g$ and $\overline{s} \geq e^{-2\sigma} s$ for some function $\sigma \in C^2 M$* and the scalar curvatures $s$ and $\bar{s}$ of $(M, g)$ and $(\overline{M}, \overline{g})$, respectively. If $\|grad\,\sigma\| \in L^1(M,g)$, then f is a homothetic mapping.

**Proof**. If $f:(M,g)\to(\overline{M},\overline{g})$ is a conformal diffeomorphism a connected complete non-compact and oriented Riemannian manifold $(M, g)$ onto another Riemannian manifold $(\overline{M},\overline{g})$ such that $\overline{g} = e^{2\sigma}g$ for some function $\sigma \in C^2 M$, then from (4.1) we obtain

$$2(n-1)\Delta\sigma = s - e^{2\sigma}\bar{s} - (n-1)(n-2)\|grad\,\sigma\|^2. \qquad (4.2)$$

Let $s \le e^{2\sigma}\bar{s}$ then (2) shows $\Delta\sigma \le 0$. It means that $\sigma$ is a superharmonic function. By the condition of our theorem, the gradient of $\sigma$ has integrable norm on $(M, g)$ and we obtain from (4.2) that $\Delta\sigma = 0$ and $\sigma$ must be harmonic (see our Lemma). Since $n \ge 3$, we see from (4.2) that $\sigma$ is constant. The proof of the theorem is complete.

Let $(M, g)$ and $(\overline{M},\overline{g})$ be Riemannian manifolds of dimension $n$ and $m$ for $n > m$. A submersion $f:(M,g)\to(\overline{M},\overline{g})$ is called a *horizontal conformal* if $f_*$ restricted to the horizontal distribution $\mathcal{H} = (Ker\,f_*)^\perp$ is conformal mapping.

Next, we consider a horizontal conformal submersion $f:(M,g)\to(\overline{M},\overline{g})$ for the case $m < n$. We note here that horizontal conformal mappings were introduced by Ishihara [25]. From the above discussion, one can conclude that the notion of horizontally conformal mappings is a generalization of concept of Riemannian submersions. In addition, we note that a natural projection onto any factor of a *double-twisted product* $(M_1 \times M_2, \lambda_1^2 g_1 + \lambda_2^2 g_2)$ of any Riemannian manifolds $(M_a, g_a)$ and smooth positive functions $\lambda_a : M_1 \times M_2 \to \mathbb{R}$ for an arbitrary $a = 1, 2$ is horizontal conformal submersion with umbilical fibres (see [20]).

Let $f:(M,g)\to(\overline{M},\overline{g})$ be a horizontal conformal submersion and $Ker\,f_*$ be an umbilical distribution then (2.2) can be rewrite in the form

$$div(\xi_V + \xi_H) = s_{mix} - \|F_H\|^2 - \frac{n-m-1}{n-m}\|\xi_V\|^2 + \frac{m-1}{m}\|\xi_H\|^2. \qquad (4.3)$$

In this case, we can formulate a corollary of Theorem 4 which generalizes our theorem on the horizontal conformal submersions of compact Riemannian manifolds with non positive mixed scalar curvature that has been proved in [26] (see also [27]).

**Corollary 5**. *Let (M, g) be an n-dimensional complete, noncompact and simply connected Riemannian manifold and $f : (M, g) \to (\overline{M}, \overline{g})$ be a horizontal conformal submersion with umbilical fibres onto another m-dimensional (m < n) Riemannian manifold $(\overline{M}, \overline{g})$. If the mean curvature vector $\xi_V$ of $\mathrm{Ker}\, f_*$ and the mean curvature vector $\xi_H$ of $(\mathrm{Ker}\, f_*)^\perp$ satisfy the condition $\|\xi_V + \xi_H\| \in L^1(M, g)$ and the mixed scalar curvature $s_{\mathrm{mix}}$ of (M, g) is nonpositive then $(\mathrm{Ker}\, f_*)^\perp$ is integrable and (M, g) is isometric to a direct product $(M_1 \times M_2,\, g_1 \oplus g_2)$ of some Riemannian manifolds $(M_1, g_2)$ and $(M_1, g_2)$ such that integral manifolds of $\mathrm{Ker}\, f_*$ and $(\mathrm{Ker}\, f_*)^\perp$ correspond to the canonical foliations of the product $M_1 \times M_2$.*

## 5. Applications to the theory of Riemannian submersions

A submersion $f : (M, g) \to (\overline{M}, \overline{g})$ is called *Riemannian submersion* if $(f_*)_x$ preserves the length of the horizontal vectors at each point $x \in M$ (see [13, p. 3]). In this case, the horizontal distribution $\mathcal{H} = (\mathrm{Ker}\, f_*)^\perp$ is totally geodesic (see [17]). In the paper [28] and in the monograph [13, p. 235] was proved the following theorem. Let $f : (M, g) \to (\overline{M}, \overline{g})$ be a Riemannian submersion with totally umbilical fibres. If (M, g) is a closed and orientable manifold with nonpositive *mixed sectional curvature* (i.e. sec (X,Y) ≤ 0 for every horizontal vector field X and for every vertical vector field Y), then all fibres are totally geodesic and horizontal distribution $\mathcal{H} = (\mathrm{Ker}\, f_*)^\perp$ is integrable, and the mixed sectional curvature is equals to zero. We present a generalization of this theorem. The following result is deduced immediately from Corollary 6.

**Corollary 6**. *Let (M, g) be an n-dimensional complete, noncompact and simply connected Riemannian manifold and $(\overline{M}, \overline{g})$ be another m-dimensional (m < n) Riemannian manifold and $f : (M, g) \to (\overline{M}, \overline{g})$ be a Riemannian submersion with totally umbilical fibres. If the mixed scalar curvature $s_{\mathrm{mix}}$ is nonpositive and the mean curvature vector $\xi_V$ of fibres satisfies the condition $\|\xi_V\| \in L^1(M, g)$, then the horizontal distribution $(\mathrm{Ker}\, f_*)^\perp$ is integrable and the Riemannian manifold (M, g) is isometric to a direct product $(M_1 \times M_2,\, g_1 \oplus g_2)$ of some Riemannian manifolds $(M_1, g_2)$ and $(M_2, g_2)$ such that the*

*integral manifolds of $\mathrm{Ker}\, f_*$ and $(\mathrm{Ker}\, f_*)^\perp$ correspond to the canonical foliations of the product $M_1 \times M_2$.*

We know from [29] that there are no Riemannian submersions from closed Riemannian manifolds with positive Ricci curvature to Riemannian manifolds with nonpositive Ricci curvature. The following statement is a direct consequence of Corollary 3 and complements the above vanishing theorem.

**Corollary 7**. *Let (M, g) be an n-dimensional complete, noncompact and simply connected Riemannian manifold and $f:(M,g)\to(\overline{M},\overline{g})$ be a Riemannian submersion onto an (n – 1)-dimensional Riemannian manifold $(\overline{M},\overline{g})$. If the vertical Ricci curvature of (M, g) is nonpositive and the mean curvature vector $\xi_V$ of fibres satisfies the condition $\|\xi_V\| \in L^1(M,g)$, then the horizontal distribution $(\mathrm{Ker}\, f_*)^\perp$ is integrable and the Riemannian manifold (M, g) is isometric to direct product $(M_1 \times M_2,\, g_1 \oplus g_2)$ of some Riemannian manifolds $(M_1, g_2)$ and $(M_2, g_2)$ such that $\dim M_1 = 1$ and the integral manifolds of $\mathrm{Ker}\, f_*$ and $(\mathrm{Ker}\, f_*)^\perp$ correspond to the canonical foliations of the product $M_1 \times M_2$.*

## 5. Applications to the theory of harmonic submersions of Riemannian manifolds

A smooth mapping $f:(M,g)\to(\overline{M},\overline{g})$ is said to be *harmonic* if $f$ provides an extremum of the energy functional $E_\Omega(f) = \int_\Omega \|f_*\|^2 d\mathrm{Vol}_g$ for each relatively closed open subset $\Omega \subset M$ with respect to the variations of $f$ that are compactly supported in $\Omega$.

If (M, g) is an n-dimensional Riemannian manifold and $f:(M,g)\to(\overline{M},\overline{g})$ is a harmonic submersion onto another m-dimensional (m < n) Riemannian manifold $(\overline{M},\overline{g})$ then each its fibre (M′, g′) is an (n – m)-dimensional closed imbedded minimal submanifold of (M, g) (see [30]). Then from the above arguments and Theorem 4 we conclude that the following corollary is true.

**Corollary 9**. *Let (M, g) be an n-dimensional complete, noncompact and simply connected Riemannian manifold and $f:(M,g)\to(\overline{M},\overline{g})$ be a harmonic submersion onto another m-dimensional (m < n) Riemannian manifold $(\overline{M},\overline{g})$ with connected fibres. If the hori-* zontal distribution $(Ker\, f_*)^\perp$ is integrable and the mixed scalar curvature $s_{\mathrm{mix}}$ is non-negative, then (M, g) is isometric to a direct product $(M_1 \times M_2,\, g_1 \oplus g_2)$ of some Riemannian manifolds $(M_1, g_2)$ and $(M_2, g_2)$ such that the integral manifolds of $Ker\, f_*$ and $(Ker\, f_*)^\perp$ correspond to the canonical foliations of the product $M_1 \times M_2$.

Using the Remark 1 we get the following

**Corollary 8**. *Let (M, g) be an n-dimensional Riemannian manifold and $f:(M,g) \to (\overline{M}, \overline{g})$ be a submersion onto another m-dimensional (m < n) Riemannian manifold $(\overline{M}, \overline{g})$ with connected fibres. If the horizontal distribution $(Ker\, f_*)^\perp$ is integrable and the mixed scalar curvature $s_{\mathrm{mix}}$ is positive, then $f:(M,g) \to (\overline{M}, \overline{g})$ is not harmonic.*

**References**

[1] Wu H.H., The Bochner technique in differential geometry, Harwood Acad. Publ., Harwood (1987).

[2] Stepanov S.E., Riemannian almost product manifolds and submersions. Journal of Mathematical Sciences, 99:6 (2000), 1788-1831.

[3] Pigola S., Rigoli M., Setti A.G., Vanishing and Finiteness Results in Geometric Analysis. A Generalization of the Bochner Technique, Birkhäuser Verlag AG, Berlin (2008).

[4] Koboyashi S., Nomizu K., Foundations of differential geometry, Volume I, Interscience Publishers, New York, 1963.

[5] Pigola S., Setti A.G., Global divergence theorems in nonlinear PDEs and geometry, Ensaios Matemáticos, 26 (2014), 1-77.

[6] Gaffney M.P., A special Stkes's theorem for complete Riemannian manifolds, Annals of Mathematics, Second Series, 60:1 (1954), 140-145.

[7] Karp L., On Stokes' theorem for noncompact manifolds, Proceedings of the American Mathematical Society, 82:3 (1981), 487-490.

[8] Caminha A., Souza P., Camargo F., Complete foliations of space forms by hypersufaces, Bull. Braz. Math. Soc., New Series, 41:3 (2010), 339-353.

[9] Caminha A., The geometry of closed conformal vector fields on Riemannian spaces, Bull. Braz. Math. Soc., New Series, 42:2 (2011), 277-300.

[10] Yau S.T., Some function-theoretic properties of complete Riemannian manifolds and their applications to geometry, Indiana Univ. Math. J., 25 (1976), 659-670.

[11] Rocamora A.H., Some geometric consequences of the Weitzenböck formula on Riemannian almost-product manifolds; weak-harmonic distributions, Illinois Journal of Mathematics, 32:4 (1988), 654-671.

[12] Reinhart B.L., Differential Geometry of Foliations, Springer Verlag, Berlin – New York, 1983.

[13] Falcitelli M., Ianus S., Pastore A.M., Riemannian submersions and related topics, Word Scientific Publishing, Singapore, 2004.

[14] Walczak P.G., An integral formula for a Riemannian manifold with two orthogonal complementary distributions, Colloquium Mathematicum, LVIII:2 (1990), 243-252.

[15] Stepanov S.E., An integral formula for a Riemannian almost-product manifold, Tensor, N. S., 55 (1994), 209-214.

[16] Stepanov S.E., Bochner's technique in the theory of Riemannian almost product structures, Mathematical notes of the Academy of Sciences of the USSR, 1990, 48:2, 778-781.

[17] Stepanov S.E., A class of Riemannian almost-product structures, Soviet Mathematics (Izv. VUZ), 33:7 (1989), 51-59.

[18] Luzynczyk M., Walczak P., New integral formula for two complementary orthogonal distributions on Riemannian manifolds, Annals of Global Analysis and Geometry, 48 (2015), 195-209.

[19] Hebda J.J., Projective maps of rank ≥ 2 are strongly projective, Differential geometry and its applications, 12 (2000), 271-280.

[20] Fernández-López M., García-Río E., Kupeli D.N., Ünal B., A curvature condition for a twisted product to be a warped product, Manuscripta Mat., 106:2 (2001), 213-217.

[21] Eisenhart L. P., Riemannian Geometry, Princeton University Press, Princeton, 1949.

[22] Giachetta G., Mangiartti I., Sardanashvily G., New Lagrangian and Hamiltonian Methods in Field Theory, Word Scientific Publishing, Singapore (1997).

[23] Stepanov S.E., On the global theory of projective mappings, Mathematical Notes, 58:1 (1995), 752-756.

[24] Stepanov S.E., Geometry of projective submersions of Riemannian manifolds, Russian Mathematics (Iz. VUZ), 43:9 (1999), 44-50.

[25] Ishihara T. A mapping of Riemannian manifolds which preserves harmonic functions, J. Math. Kyoto Univ., 19 (1979), 215-229.

[26] Stepanov S.E., Weyl submersions, Russian Mathematics (Izvestiya VUZ, Matematika), 36:5 (1992), 87-89.

[27] Zawadzki T., Existence conditions for conformal submersions with totally umbilical fibers, Differential Geometry and its Applications, 35 (2014), 69-85.

[28] Bădiţoiu G., Ianuş S., Semi-Riemannian submersions with totally umbilical fibres, Rendiconti del Circolo Matematico di Palermo, 51:2 (2002), 249-276.

[29] Pro C., Wilhelm F., Riemannian submersions need not preserve positive Ricci curvature, Proc. Amer. Math. Soc., 142:7 (2014), 2529-2535.

[30] Stepanov S.E., $O(n) \times O(m-n)$-structures on $m$-dimensional manifolds, and submersions of Riemannian manifolds, St. Petersburg Math. J., 7:6, 1005-1016 (1996).